\documentclass[11pt]{article}

\usepackage{amssymb,amsmath}
\usepackage[T1]{fontenc}

\usepackage[english]{babel}

\newcommand{\Pq}[1]{\left[ #1 \right] }

\renewcommand{\sinh}{\mathrm{Sinh}}
\renewcommand{\tanh}{\mathrm{Tanh}}

\newcommand{\bproof}{\begin{proof}}
\newcommand{\eproof}{\end{proof}}

\newcommand{\bcas}{\begin{cases}}
\newcommand{\ecas}{\end{cases}}

\renewcommand{\b}[1]{{\bf #1}}

\newcommand{\SR}{{sr}}

\renewcommand{\b}[1]{{\bf #1}}

\newcommand{\Dh}{\Delta_\SR}

\newcommand{\funz}[5]{#1 : \begin{tabular}{ccl}
 #2 &$\rightarrow$& #3 \\
 #4 & $\mapsto$& #5   \end{tabular}}

\newcommand{\mmfunz}[5]{\begin{center}
\funz{$\displaystyle{#1}$}{$\displaystyle{#2}$}{$\displaystyle{#3}$}{$\displaystyle{#4}$}{$\displaystyle{#5}$}
\end{center}}

\newcommand{\llabel}[1]{{\label{#1}}}

\newcommand{\ffoot}[1]{}
\newcommand{\fffoot}[1]{}

\renewcommand{\r}[1]{(\ref{#1})}

\newcommand{\bi}{\begin{itemize}}
\newcommand{\ei}{\end{itemize}}

\newcommand{\be}[1][]{\begin{enumerate}[#1]}
\newcommand{\ee}{\end{enumerate}}

\newcommand{\bd}{\begin{description}}
\newcommand{\ed}{\end{description}}
\renewcommand{\i}{\item}

\newcommand{\bqn}{\begin{eqnarray}}
\newcommand{\eqn}{\end{eqnarray}}
\newcommand{\eqnn}{\nonumber\end{eqnarray}}
\newcommand{\eqnl}[1]{\llabel{#1}\end{eqnarray}}

\newcommand{\nn}{\nonumber\\}
\newcommand{\ba}[1]{\begin{array}{#1}}
\newcommand{\ea}{\end{array}}

\newcommand{\R}{\mathbb{R}}
\newcommand{\C}{\mathbb{C}}

\newtheorem{Theorem}{\bf Theorem}[section]
\newtheorem{lemma}[Theorem]{\bf Lemma}
\newtheorem{corollary}[Theorem]{\bf Corollary}
\newtheorem{definition}[Theorem]{\bf Definition}
\newtheorem{proposition}[Theorem]{\bf Proposition}
\newtheorem{example}[Theorem]{\bf Example}
\newtheorem{remark}[Theorem]{\bf 
{{\sl Remark}}}
\newenvironment{proof}{{\it Proof.}~~}{\hfill$\square$}

\newcommand{\bt}{\begin{Theorem}}
\newcommand{\et}{\end{Theorem}}
\newcommand{\bl}{\begin{lemma}}
\newcommand{\el}{\end{lemma}}
\newcommand{\bp}{\begin{proposition}}
\newcommand{\ep}{\end{proposition}}
\newcommand{\bc}{\begin{corollary}}
\newcommand{\ec}{\end{corollary}}
\newcommand{\bdeff}{\begin{definition}}
\newcommand{\edeff}{\end{definition}}
\newcommand{\brem}{\begin{remark}\rm}
\newcommand{\erem}{\end{remark}}
\newcommand{\bex}{\begin{example}\rm}
\newcommand{\eex}{\hfill$\triangledown$\end{example}}

\newcommand{\lam}{\lambda}

\newcommand{\Pt}[1]{\left( #1 \right)}
\newcommand{\Pg}[1]{\left\{ #1 \right\}}

\newcommand{\auth}[1]{{\sc #1}}
\newcommand{\tit}[1]{{\rm #1}}
\newcommand{\titl}[1]{{\it #1}}
\newcommand{\jou}[1]{{\it #1}}

\newcommand{\pp}[1]{pp.~#1}

\newcommand{\G}{G}
\renewcommand{\L}{\mathfrak{L}}

\newcommand{\DD}{\mathcal{C}^\infty_c(\G)}
\renewcommand{\div}{\mathrm{div}}
\newcommand{\dil}[1][r]{\Gamma_{#1}} 
\newcommand{\bassor}{{|_{\Gamma_{r}(x)}}}

\title{Large time behavior for the heat equation on Carnot groups}

\author{Francesco Rossi}

\begin{document}
\maketitle

\begin{abstract}
We first generalize a decomposition of functions on Carnot groups as linear combinations of the Dirac delta and some of its derivatives, where the weights are the moments of the function.

\noindent We then use the decomposition to describe the large time behavior of solutions of the hypoelliptic heat equation on Carnot groups. The solution is decomposed as a weighted sum of the hypoelliptic fundamental kernel and its derivatives, the coefficients being the moments of the initial datum.
\end{abstract}

\section{Introduction}

The study of Carnot groups and of PDEs on these spaces have drawn an increasing attention for several reasons. First, Carnot groups are topologically extremely simple, since they are isomorphic to $\R^N$, but in the same time their metric structure is different, since they are naturally endowed with Carnot--Caratheodory metrics. Second, several tangent spaces to Carnot--Caratheodory metric spaces are diffeomorphic to Carnot groups. Third, and more connected to PDEs, Carnot groups are the easiest examples of spaces where hypoelliptic equations are naturally defined (see \cite{lanconelli-book} and references therein). Thus, they are the most natural examples in which one can study relations between the solutions of the hypoelliptic PDEs and the Carnot--Caratheodory metrics.

With this goal, we first generalize a decomposition of functions first stated in \cite{moments} for functions on $\R^N$. In that article, tha authors find conditions to write functions as
\bqn
f=a_0\delta_0+\sum_{i=1}^N a_i \partial_i \delta_0+\sum_{i,j=1}^N a_{ij} \partial_{ij} \delta_0+\ldots,
\eqnl{e-dec}
where $\delta_0$ is the Dirac delta at 0 and $\delta_i,\delta_{ij},\ldots$ are standard derivatives in $\R^N$. The coefficients $a_0, a_i,a_{ij},\ldots$ are called the moments. Our goal is to extend such kind of formula to Carnot groups. In this case, we have to replace $\R^N$ derivatives with derivatives that are ``adapted'' to the Carnot group. The natural choice is to use derivative operators that are invariant with respect to the operation of the Carnot group. We will recall all the necessary definitions and result for Carnot groups in Section \ref{s-carnot}.\\

One of the interest of \r{e-dec} is that it is useful to study the heat equation on Carnot groups. In this case, the standard Laplacian $\Delta=\sum_{i=1}^N {\delta_i^2}$ is replaced by the sub-Riemannian Laplacian $\Dh$, that is an hypoelliptic operator. Also in this case, all necessary definitions and result are recalled in Section \ref{s-carnot}. In particular, we study the initial value problem
\bqn
\label{e-cauchy}
\begin{cases}
\partial_t f=\Dh f \\
f(t=0)=f_0,
\end{cases}
\eqn
whose solution is given by $f(t)=f_0\ast P_t$, where the convolution $\ast$ is intrinsically defined on the Carnot group $\G$ by its Lie group structure (see Definition \ref{def-conv}) and $P_t$ is the fundamental solution of the hypoelliptic heat equation with initial datum $\delta_0$.  More precisely, we study the large time behavior of solutions of \r{e-cauchy} and present a complete asymptotic expansion using the fundamental solution and its derivatives as basis functions and the moments of the initial datum as coefficients.  We show that, due to anisotropicity, only a subset of the moments needs to be used.

The structure of the paper is the following. In Section \ref{s-carnot}, we recall the definition of Carnot groups, of the metrics defined on them, the convolution and some useful properties of the fundamental solutions for the heat equation. Section \ref{s-moments} is devoted to the main result of this article, the decomposition of a function with respect to its moments. Section \ref{s-pde} is devoted to the study of the initial value problem \r{e-cauchy}, with a description of the large time behaviour in terms of the moments of the initial datum. We finally apply these results to the case of the Heisenberg group.

\section{Carnot groups and hypoelliptic heat equations}
\label{s-carnot}

\noindent We define here Carnot groups and recall some of their fundamental properties. For more details, see \cite{lanconelli-book,folland}.

Let $\L$ be a nilpotent Lie algebra that, by definition, is a vector space, thus isomorphic to $\R^N$ for some $N \ge 1$. A Lie algebra is endowed with a Lie bracket product $[\cdot,\cdot]$. One can always assume that the Lie algebra is an algebra of matrices (Ado theorem, see e.g. \cite{barut}), and that the Lie bracket is $[A,B]=AB-BA$. We recall that the Lie algebra is nilpotent where there exists a $n>0$ such that $[A_1,[A_2,[A_3,\ldots [A_{n-1},A_n]]..]]=0$ for any choice of $A_1,\ldots,A_n\in\L$. In terms of matrices, a nilpotent algebra can be always identified with a subalgebra of the algebra of upper triangular matrices of a given dimension with zero diagonal.

Let $\G$ be the exponential of $\L$. Again, in terms of matrices, one can consider the exponential map $\exp:\L\rightarrow \G$ to be the standard matrix exponential. $G$ is a group with respect to the standard matrix product, and it is indeed a Lie group. For more details, see e.g. \cite{barut}. Since $\L$ is nilpotent, the exponential map $\exp:\L\rightarrow \G$ is a diffeomorphism. We fix a basis\footnote{The notation $X^i$ with $i\in\Pg{1,\ldots,N}$ denotes, all along the paper, a vector of the basis of the Lie algebra.} $\Pg{X^1,\ldots, X^N}$ on $\L$. All along the paper, the notation $x=(x_1,\ldots,x_N)$ can denote both $x=\sum x_i X^i\in\L$ or $x=\exp\Pt{\sum x_i X^i}\in \G$, depending on the context. Also recall that the Lie algebra of a Lie group is identified with the space of invariant differential operators on the Lie group, i.e. of vector fields $X$ such that $X(g)=(L_g)_*\,  X(e)$, where $L_g$ is the left translation $h\to gh$ and $(L_g)_*$ is its differential. See more details in \cite{barut}.

We consider the Lebesgue measure $dx=dX^1\ldots dX^N$. This gives to $\G$ a (bi-invariant) Haar measure (the lift of the Lebesgue measure on $\L$) on $\G$, that we still denote with $dx$. In the following, we deal with $L^p$ spaces on $\G$ with respect to this measure. We recall that the Haar measure is invariant in the sense that it satisfies
$\int f(x)\, dx=\int f(y \cdot x)\, dx$ for any $y\in \G$.

\bex It is interesting to recall that all the definitions given above are generalizations of the standard properties of the Euclidean space, i.e. $\R^n$ endowed with the vector sum $(x^1,\ldots,x^N)+(y^1,\ldots,y^N)=(x^1+y^1,\ldots,x^N+y^N)$. Indeed, $\R^n$ itself is a Lie algebra with trivial brackets $[A,B]=0$. It is realized as a matrix algebra by writing each element $(x^1,\ldots,x^N)$ as the matrix of dimension $N+1$ 
$$A=\Pt{\ba{cccc} 0&\ldots& 0 &x^1\\
0&\ldots& 0 &x^2\\
\vdots&\ldots&\vdots&\vdots\\
0&\ldots& 0 &x^N\\
0&\ldots& 0 &0
\ea}.$$
The exponential $\G$ satisfies $e^{A}=I+A$, and indeed in this case we have $e^{A}e^{B}=e^{A+B}$, thus $\G$ is isomorphic to the Euclidean group. The Haar measure is the standard Lebesgue measure on $\R^N$, and the invariance of the Haar measure is indeed the invariance of the Lebesgue measure with respect to translations.
\eex

\bex \label{ex-1}The simplest (nontrivial) example of nilpotent Lie algebra and Lie group is given by the \b{Heisenberg group} $H_2$, see e.g. \cite[Chap. 1]{lanconelli-book}.  Choose the Lie algebra $\L$ as the vector space generated by $\Pg{l^1,l^2,l^3}$ where
\bqn
l^1=\Pt{\ba{ccc}
0&1&0\\
0&0&0\\
0&0&0
\ea}\quad
l^2=\Pt{\ba{ccc}
0&0&0\\
0&0&1\\
0&0&0
\ea}\quad
l^3=\Pt{\ba{ccc}
0&0&1\\
0&0&0\\
0&0&0
\ea}.
\eqnl{e-H2-alg}
They satisfy the following commutation rules: $[l^1,l^2]=l^3$, $[l^1,l^3]=[l^2,l^3]=0$. Hence, $\L$ is a nilpotent Lie algebra. Computing the matrix exponential, one can write the Heisenberg group as the 3D group of matrices
$$H_2=\Pg{\Pt{\ba{ccc}
1&x&z+\frac12 x y\\
0&1&y\\
0&0&1
\ea}\ |\ x,y,z\in\R}$$
endowed with the standard matrix product. Left-invariant vector fields are of the form $X(g)=g l$, where $g\in H_2$ and $l\in \L$.

One can write $H_2$ as $\R^3,$ via the identification
$$(x,y,z)\sim\Pt{\ba{ccc}
1&x&z+\frac12 x y\\
0&1&y\\
0&0&1
\ea}.$$ In this case, the group law is $$(x_1,y_1,z_1)\cdot(x_2,y_2,z_2)=\Pt{x_1+x_2,y_1+y_2,z_1+z_2+\frac12\Pt{x_1y_2-x_2y_1}}.$$ The corresponding left-invariant vector fields are generated by 
$${X^1}=\partial_x-\frac y2 \partial_z,\qquad{X^2}=\partial_y+\frac x2 \partial_z,\qquad X^3=\partial_z.$$
\eex

A family of dilations on $\L$ is a one-parameter family $\Pg{\dil\ |\ 0<r<\infty}$ of automorphisms of $\L$ of the form $\dil=\exp(A \log r)$, where $A$ is a diagonalizable linear transformation of $\L$ with positive eigenvalues.  Dilations of $\L$ lift to dilations of $\G$, that we still denote by $\dil$. The constant $Q=\mathrm{Tr}(A)$ is called {\it homogeneous dimension} of $\G$. A measurable function $f$ on $\G$ is called {\it homogeneous of degree $\lam$}  if $f\circ \dil=r^\lam f$ for all $r>0$. A differential operator $D$ is called {\it homogeneous of degree $\lam$} if $D(u\circ \dil)=r^\lam (Du)\circ \dil$ for all $u\in\DD$ and $r>0$. We have $d(\dil x) = r^Q dx.$  For proofs and further results, see \cite{folland}.

We focus on a class of nilpotent Lie algebras, the \b{stratified Lie algebras}. A stratification of a Lie algebra $\L$ is a decomposition $\L=\oplus_{i=1}^s \L_i$ satisfying $\Pq{\L_1,\L_i}=\L_{i+1}$ for $1\leq i < s$ and $\Pq{\L_1,\L_s}=0$. If $\L$ is stratified, then it admits a family of dilations
$$\dil(X_1+\ldots+X_s)=r X_1+r^2 X_2+\ldots+r^s X_s,$$
where\footnote{The notation $X_i$ with $i\in\Pg{1,\ldots,s}$ denotes, all along the paper, a generic vector such that $X_i\in\L_i$.} $X_i\in \L_i$. A \b{stratified group} $\G$ is the exponential of a stratified Lie algebra. In this case we have $Q=\sum_{i=1}^s i\, \mathrm{dim} (\L_i)$. Given a vector field $X$ and considering it as a differential operator, it is homogeneous of degree $i$ if and only if $X\in \L_i$.  From now on, we denote $s$ the step of nilpotency of $\G$ and $w^i$ the degree of homogeneity of $X^i$. We also denote $n={\dim(\L_1)}$.

\bex Once again, the Euclidean group is a particular case of the stratified algebras, in which the dilation is the standard ``isotropic'' dilations $\dil(x^1,\ldots,x^N)=(rx^1,\ldots,rx^N)$. The homogeneous dimension coincides with $N$. Elements of $\L$ are considered as derivatives simply via the idenfitication $x^i\to \partial_i$, and they are indeed homogenoeus of order 1 with respect to the dilation.

The basic idea of Carnot groups is, more in general, that dilations can have different degrees of homogeneity with respect to different axes, that corresponds to ``anisotropic'' dilations.
\eex

\bex \label{ex-2} The Heisenberg group is a stratified Lie group, choosing $\L_1=\mathrm{span}\Pg{X^1,X^2}$ and $\L_2=\mathrm{span}\Pg{X^3}$. Its homogeneous dimension is $4$.
\eex

Given a stratified Lie group $\G$, we define a quadratic form $<.,.>_S$ on $\L_1$, a structure that is called a \b{Carnot group}. Carnot groups are naturally endowed with a metric structure, the so called Carnot--Caratheodory metric. We don't use this metric in this article. For further references, see e.g. \cite{lanconelli-book}.

In the following we use another norm on a Carnot group. Define a quadratic form $<.,.>_R$ on the whole Lie algebra $\L$ such that it coincides with $<.,.>_S$ on $\L_1$ and such that the $\L_i$ are orthogonal with respect to this quadratic form. We introduce the corresponding norm $|v|_R:=<v,v>_R$ on $\L$ and define the following norm on $\G$, that is homogeneous of degree 1.
\bqn
\|\exp\Pt{X_1+X_2+\ldots+X_s}\|_R&:=&\Pt{\sum_{i=1}^s |X_i|_R^{\frac{2(s!)}i}}^{\frac{1}{2(s!)}}
\eqnn
\bex \label{ex-3}
The Heisenberg group is a Carnot group when endowed with a quadratic form on $\L_1$. The standard choice is given by choosing $X^1,X^2$ orthonormal. We also endow the whole $\L$ with a quadratic form, by choosing $X^1,X^2,X^3$ orthonormal. Thus $\|(x,y,z)\|_R=\sqrt[4]{|x|^{4}+|y|^{4}+|z|^2}.$
\eex

It is interesting to remark that the quadratic form $<.,.>_R$ is not intrinsic on the Carnot group $(\G,<.,.>_S)$. Nevertheless, all the possible norms are equivalent, i.e. given two quadratic forms $<.,.>_{R_1}$ and $<.,.>_{R_2}$ on $\L$, there exist $c,C>0$ such that $c<v,v>_{R_1}\leq<v,v>_{R_2}\leq C<v,v>_{R_1}$. We assume that the vectors $\Pg{X^1,\ldots,X^n}$ are orthonormal with respect to $<.,.>_S$ and that $\Pg{X^1,\ldots,X^N}$ are orthonormal with respect to $<.,.>_R$.





\subsection{Hypoelliptic heat equations on Carnot groups}

As already stated, one of the main features of Carnot groups is that they are naturally endowed with a Laplacian operator. Given $\Pg{X^1,\ldots,X^n}$ an orthonormal basis of $\L_1$ with respect to $<.,.>_S$, the intrinsic Laplacian for $\G$ is the sum of squares
\bqn
\mbox{$\Dh:=\sum_{i=1}^n {X_i}^2.$}
\eqnn
Here, we are interested in the \b{hypoelliptic heat operator} $\partial_t-\Dh$.

The Laplacian and the heat operator are intrinsic in the sense that they are intrinsically defined from the structure of the Carnot group. In particular, one can find connections between the solution of the hypoelliptic heat equation and the intrinsic (sub-Riemannian) distance defined on Carnot group. This interesting topic is outside the goals of this paper: more informations and interesting estimates can be found in \cite{nostro-kern,cyg,eldredge,varop}. Other results for the hypoelliptic heat operator defined on other spaces than Carnot groups can be found in \cite{alex,ostellari}.

We recall that both the Laplacian and the heat operator are hypoelliptic, according to the following definition (see e.g. \cite{hormander}).
\bdeff
An operator $L$ is hypoelliptic when, given $U\subset\R^N$ and $\phi:U\to\C$ such that $L\phi\in C^\infty$, then $\phi$ is $C^\infty$.
\edeff

We are interested in the \b{hypoelliptic heat equation} on Carnot groups. We recall some properties.
\bt[\cite{folland} (2.1),(2.11),(3.1)] Let $\G$ be a Carnot group of homogeneous dimension $Q>2$. Then the {\it heat operator} $\partial_t-\Dh$ is homogeneous of degree 2. There is a unique heat kernel \mmfunz{P}{\G\times(0,+\infty)}{\R}{(x,t)}{P_t(x)} of type 2 such that $P_0(.)=\delta_0$. It satisfies $P_t(x)\geq0$, $\int_\G P_t(x) dx=1$ for all $t$ and $P_{r^2 t}(\dil(x))=r^{-Q} P_t(x)$.
\et

For the following, we need estimates about the large time behavior of the heat kernel $P_t$ and of its derivatives (of any order) with respect to left-invariant vector fields. These results are proved in \cite{folland}.
\bc[\cite{folland} (3.4)-(3.6)] The kernel $P_t(x)$ is $C^\infty$ on $\G\times(0,+\infty)$. For each $t_0>0$ and positive integer $k$ there is a constant $C_{t_0,k}$ such that, for all $\|x\|_R\geq 1$ and $t\leq t_0$, we have
\bqn \label{p-decad}
|P_t(x)|\leq C_{t_0,k} \|x\|_R^{-k}.
\eqn

Let $D$ be any left-invariant differential operator on $\G$. Then the same estimate \r{p-decad} holds for $DP_t$.
\ec
\noindent We now give estimates for the $L^p$ norms of the heat kernel and its derivatives.
\bl \label{l-P1} Let $\G$ be a Carnot group of homogeneous dimension $Q>2$, $P_t$ the corresponding heat kernel, with $t>0$. Given $1 \leq q \leq \infty$, we have $\|P_1\|_q<\infty$ and $\| P_t\|_q= t^{-Q/(2q^*)} \| P_1\|_q$, with $\frac1q+\frac1{q^*}=1$.

Given an homogeneous left-invariant differential operator $D$ of degree $\lam$, we have $\|DP_1\|_q<\infty$ and $\| DP_t\|_q= t^{-Q/(2q^*)-\lam/2} \| DP_1\|_q$.
\el
\bproof
For $1<q<\infty$ we have $ \| P_t\|^q_q=\int_\G |P_t(x)|^q dx=\int_\G r^{Qq} |P_{r^2t}(\dil(x))|^q dx=r^{Qq} \int_\G |P_{r^2t}(y)|^q r^{-Q}dy=r^{Q(q-1)}\|P_{r^2t}\|^q_q.$
The proof is given by choosing $r=t^{-\frac12}$. Remark that $\|P_1\|_q$ is finite due to \r{p-decad}. The proof for $q=1$ and $q=\infty$ is equivalent. The proofs for $DP_t$ are identical, recalling $DP_{r^2t}(\dil x)=r^{-Q-\lam}DP_t(t,x)$.
\eproof

\subsection{Solutions of the initial value problem}

Given the fundamental solution $P_t$ for the hypoelliptic heat equation, we want to find the solution of the initial value problem \r{e-cauchy} for a given $f_0$. The solution is given by the \b{convolution} $f_t=f_0\ast P_t$, with the convolution defined on Lie groups as follows.
\bdeff
\label{def-conv}
Let $(G,\cdot)$ be a Lie group and $dg$ its left-invariant Haar measure. The convolution is
\bqn
(f_1\ast f_2)(g)=\int_G  f_1(h)f_2(h^{-1}\cdot g)\, dh.
\eqnn
\edeff
Remark that the product and the inverse are given by the group operation on $\G$. The convolution on Lie groups satisfies an invariance property with respect to left-invariant differential operators. More precisely, given a left-invariant vector field $X$, we have 
\bqn
X (f_1)\ast f_2=X (f_1\ast f_2)=f_1\ast X (f_2).
\label{e-convder}
\eqn
\brem The standard definition of convolution on $\R^N$ is given by considering it as the Euclidean space. In this case, left-invariant vector fields are differential operators with constant coefficient.
\erem
A useful property is the Young's inequality
\bt[\cite{hewitt} (20.14-18)]
Let $f_1\in L^1(\G)$ and $f_2 \in L^p(\G)$. Then $f_1\ast f_2 \in L^p(\G)$ and 
\bqn
\|f_1\ast f_2\|_p\leq \|f_1\|_1 \|f_2\|_p.
\eqnn
\et

\section{Functions decomposition on Carnot groups}
\label{s-moments}

In this section, we generalize to Carnot groups a decomposition of functions based on Dirac delta and its derivatives, first stated in \cite{moments} for functions defined on $\R^N$. We first recall that, given $F=(F_1,\ldots,F_N)$, the intrinsic divergence\footnote{Also in this case, the divergence is intrinsic in the sense that is intrinsically defined by the structure of the Carnot group. For more details, see e.g \cite{nostro-kern}.}  is $\mathrm{div} F=\sum_{i=1}^N X^i(F_i)$.

We give a first decomposition of functions, based on Dirac delta only. The result is not interesting in itself, since it is weaker than \cite[Thm. 1]{moments}, but its proof gives an example of the method used in the following to deal with Carnot groups.

\bp \label{t:ez1} Let $f$ be a function of $G$ and $p\in[1,\infty]$.
\bi
\i If $ p < Q/ (Q-1)$, $f\in L^1(\G)$ and $|x|_Rf\in L^p(G)$, then it exists $F\in \Pt{L^p(G)}^N$ such that $f=\Pt{\int f} \delta_0 + \div F$.
\i If $p > Q/ (Q-s)$ and $|x|_Rf\in L^p(G)$, then it exists $F\in \Pt{L^p(G)}^N$ such that $f=\div F$.
\ei
We also have that $\|F_j\|_p\leq C_{p,N} \| |x|_R f\|_p$ for each $j=1,\ldots,N$.
\ep
\bproof The proof is a slight generalization of the proof given in \cite{moments}. Fix $x=(x_1,\ldots,x_N)\in\G$ and take a test function $\phi\in\mathfrak{D}(\G)$. For the first case, we apply the fundamental theorem of integral calculus for the function of $r\in[0,1]$ defined by $\tilde\phi(r):=\phi(\dil(x))$ finding 
\bqn
\phi(x)-\phi(0)&=&\int_0^1 \frac{\partial \tilde\phi}{\partial r}\,dr=\int_0^1 \sum_{i=1}^N X^i(\phi)_\bassor \frac{\partial (r^{w^i}x_i)}{\partial r}\,dr=\nn
&=&\int_0^1 \sum_{i=1}^N w^i r^{w^i-1} x_i X^i(\phi)_\bassor\, dr.
\eqnn
We now apply $f-(\int f)\delta_0$ to $\phi(x)$ and have
\bqn
\int_\G\Pt{f-(\mbox{$\int$} f)\delta_0}\,\phi\, dx&=&\int_\G f(x)\int_0^1 \sum_{i=1}^N  w^i x_i r^{w^i-1} X^i(\phi)_\bassor \, dr dx=\nn
&=&-\int_\G\phi(x)\sum_{i=1}^N X^i\Pt{ w^i x_i\int_0^1 \frac{f(\Gamma_{r^{-1}}(x))}{r^{Q+1} }\,dr}\, dx.
\eqnn
The last identity is given by integration by parts, since the $X^i$ are skew-adjoint with respect to $dx=dX^i\ldots dX^N$, see e.g. \cite[Pr.20]{nostro-kern}. We also apply change of variables $\dil(x)\mapsto x'$. Choose $F_i(x):= w^i x_i \int_0^1 \frac{f(\Gamma_{r^{-1}}(x))}{r^{Q+1} }\,dr$. In this case we have 
\bqn
\|F_i\|_p^p = |w^i|^p\cdot\| |x|_R\, f(x)\|_p^p\cdot |\int_0^1 r^{\frac{Q}p +w^i -Q-1}dr|^p.
\eqnl{e-defF} It is finite when $p< \frac{Q}{Q-w^i}$. Since $w^i\geq 1$, it is finite when $p<\frac{Q}{Q-1}$.

The other case is similar, given the formula 
$$\phi(x)=-\int_0^1 \sum_{i=1}^N w^i r^{-1-w^i} x_i X^i(\phi(\Gamma_{r^{-1}}(x))\,dr,$$
from which we find 
\bqn
F_i(x):= w^i x_i \int_0^1 r^{Q-1} f(\dil(x))\, dr.
\eqnl{e-defF2}
In this case, we have $$\|F_i\|_p^p = |w^i|^p\cdot\| |x|_R\, f(x)\|_p^p\cdot |\int_0^1 r^{-\frac{Q}p +Q-1-w^i}dr|^p.$$ It is finite when $p>\frac{Q}{Q-w_i}$. Since $w_i\leq s$, it is finite when $p>\frac{Q}{Q-s}$.
\eproof

\brem \label{rem-RN} One can check that all the formulas reduce to the results of \cite{moments} when $G=(R^N,+)$, recalling that $s=1$, $w_i\equiv 1$, $\dil(x)=rx$ and $Q=N$. The main difference with \cite{moments} is the set of $L^p$ spaces in which none of the decompositions holds, that is $p\in\Pq{\frac{Q}{Q-1},\frac{Q}{Q-s}}$.
\erem

\brem \label{rem-otherF} One can ask if a similar decomposition can hold for $p\in\Pq{\frac{Q}{Q-1},\frac{Q}{Q-s}}$. It is clearly possible, since one can take results of \cite[Thm. 1]{moments}, provided $p\neq \frac{N}{N-1}$. Remark that $\frac{N}{N-1}\in\Pq{\frac{Q}{Q-1},\frac{Q}{Q-s}}$. 



If one looks for a formula $f=\Pt{\int f} \delta_0 + g$ with $g$ given by a differential operator that is not the divergence, it is certainly possible to find more complicated decomposition. We use this idea in the following, that gives operators in which vector fields $X^i$ play different roles depending on their weight $w_i$.
\erem

As stated above, a higher order decomposition can be defined using the derivatives of the Dirac delta as basis functions, following \cite[Th. 2]{moments}. In this case, it is necessary to compute high order derivatives of $\tilde\phi(s)$. For this reason, we only state the first order decomposition (its meaning will be clear in the following). For higher order decomposition one has to apply the same technique with longer computations.

\bt \label{t:ez2} 
Let $f$ be a function of $G$ and $p\in[1,\infty]$.
\bi
\i If $p<\frac{Q}{Q-1}$, $f\in L^1(\G,1+|x|_R)$ and $|x|_Rf,|x|_R^2 f\in L^p(\G)$, then
\bqn
f=(\int_G f)\delta_0 -\sum_{i=1}^n(\int f x_i) (X_i \delta_0)+\sum_{i,j=1}^N X^i X^j(F_{ij})+\sum_{i=n+1}^N X^i(F_i).
\eqnn
\i If $p>\frac{Q}{Q-2s}$ and $|x|_Rf,|x|_R^2 f\in L^p(\G)$, then
\bqn
f=\sum_{i,j=1}^N X^i X^j(F_{ij})+\sum_{i=n+1}^N X^i(F_i).
\eqnn
\ei
In both cases, we have $\|F_{ij}\|_p\leq C \| |x|_R^2 f\|_p$ and $\|F_i\|_p\leq C \| |x|_R f\|_p$.
\et
\bproof
Fix $x=(x_1,\ldots,x_N)\in\G$. Given a test function $\phi\in\mathfrak{D}(\G)$, consider the function of $r\in[0,1]$ defined by $\tilde\phi(r):=\phi(\dil(x))$. The Taylor polynomial $P_k\tilde\phi(r)$ of order $k$ near $0$ for $\tilde\phi$ satisfies
\bqn
\tilde\phi(1)-P_k\tilde\phi(1)&=&\int_0^1 \frac{\tilde\phi^{(k+1)}(r)}{k!}(1-r)^k\,dr.
\eqnn
We also have
\bqn
\tilde\phi(1)&=&(-1)^{k+1}\int_1^\infty \frac{\tilde\phi^{(k+1)}(r)}{k!}(r-1)^k\,dr=\nn
&=&(-1)^{k+1} \int_0^1 \frac{\tilde\phi^{(k+1)}\Pt{\frac{1}{t}}}{k!}(1-t)^k\,\frac{dt}{t^{k+2}},
\eqnl{e-taylor-inf}
that is a particular case\footnote{Remark a misprint in \cite[eq. (6)]{moments}, where $(-1)^k$ must be replaced with $(-1)^{k+1}$.} of
\bqn
\tilde\phi(r)&=&(-1)^{k+1} \int_0^1 \frac{\tilde\phi^{(k+1)}\Pt{\frac{r}{t}}}{k!}\Pt{\frac{r}{t}}^{k+1}(1-t)^k\,\frac{dt}{t}.
\eqnn
We merge the two in this more general formula
\bqn
&&\tilde\phi(1)-P_{k-j}\tilde\phi(1)=\int_0^1\frac{\tilde{\phi}^{(k-j+1)}(r)}{(k-j)!} (1-r)^{k-j}\,dr=\nn
&&=(-1)^j\int_0^1 dr \int_0^1 \frac{dt}{t} \Pt{\frac{r}{t}}^j (1-t)^{j-1} (1-r)^{k-j} 
\frac{\tilde\phi^{(k+1)}\Pt{\frac{r}{t}}}{(k-j)! (j-1)!}.
\eqnl{e-taylor-kj}

We now restrict ourselves to the case $k=1$. We compute $\tilde\phi'$ and $\tilde\phi''$:
\bqn
\tilde\phi'(r)&=&\sum_{i=1}^N X^i(\phi)_\bassor w^i r^{w^i-1} x_i,\nn
\tilde\phi''(r)&=&\sum_{i,j=1}^N X^iX^j(\phi)_\bassor w^i w^j r^{w^i+w^j-2} x_i x_j +\nn
&&\hspace{1cm}+\sum_{i=n+1}^N X^i(\phi)_\bassor w^i(w^i-1) r^{w^i-2} x_i.
\eqnl{e-der2}

Remark that the last term in the r.h.s. of \r{e-der2} contains first order derivatives with respect to vector fields whose degree of homogeneity is greater than or equal to 2. We compute the Taylor formulas for $\tilde \phi$:
\bqn
P_0 \tilde\phi(1)&=&\phi(0),\qquad P_1\tilde\phi(1)=\phi(0)+\sum_{i=1}^n X^i(\phi)_{|_0} x_i.
\eqnn
Remark that $P_1\tilde\phi$ contains derivatives in $\L_1$ only.

We now apply $f-\Pt{(\int_G f) \delta_0-\sum_{i=1}^n(\int f x_i) (X_i \delta_0)}$ to $\phi$, that gives
\bqn
&&\int\Pt{f(x)\phi(x)-  f(x) \phi(0) -\sum_{i=1}^n f(x) x_i X^i(\phi)_{|_0}}\,dx=\nn
&=&\int f(x) \int_0^1 \frac{\tilde\phi''(r)}{k!}(1-r)^k\,dr\,dx=\nn
&=&\int \frac{f(x)}{k!} \Pt{\int_0^1\sum_{i,j=1}^N X^i X^j(\phi)_\bassor w^i w^j r^{w^i+w^j-2} x_i x_j +\right.\nn
&&\hspace{3mm}\left. +\sum_{i=n+1}^N X^i(\phi)_\bassor w^i(w^i-1) r^{w^i-2} x_i} (1-r)\,dr\,dx=\nn
&=&\int_G\phi(x)\Pt{\sum_{i,j=1}^N \int_0^1 w^i w^j X^j X^i(f(\Gamma_{r^{-1}}(x))x_i x_j) \frac{1-r}{r^{Q+2}}+\right.\nn
&&\left.\hspace{3mm} - \sum_{i=n+1}^N \int_0^1 w^i (w^i-1) X^i(f(\Gamma_{r^{-1}}(x))x_i) \frac{1-r}{r^{Q+2}}} dr\,dx.
\eqnn
Choose $F_{ij}=F_{ji}:=\int_0^1 w^i w^j f(\Gamma_{r^{-1}}(x))x_i x_j \frac{1-r}{r^{Q+2}}\,dr$, that satisfies
\bqn
\|F_{ij}\|_p^p= |w^i w^j|^p \int |f(x) x_i x_j|^p \,dx |\int_0^1 r^{w^i+w^j-2-Q+\frac{Q}{p}}(1-r)\,dr|^p,
\eqnn
that is finite when $f(x) x_i x_j\in L^p(\G)$ and $w^i+w^j-2-Q+\frac{Q}{p}>-1$. Since $w_i\geq 1$, it is finite when $p < \frac{Q}{Q-1}$. Similarly, for $i\geq n+1$ choose $F_i(x):=\int_0^1 w^i (w^i-1) f(\Gamma_{r^{-1}}(x))x_i \frac{1-r}{r^{Q+2}}\,dr$, that satisfies
\bqn
\|F_i\|_p^p= |w^i (w^i-1)|^p \int |f(x) x_i|^p \,dx \mid\int_0^1 r^{w^i-2-Q-\frac{Q}{p}}(1-r)\,dr\mid^p.
\eqnn
It is finite when $f(x) x_i \in L^p(\G)$ and $w^i-2-Q+\frac{Q}{p}>-1$. Since $w^i \geq 2$, it is finite when $p < \frac{Q}{Q-1}$.

The second case is similar, using \r{e-taylor-inf} and \r{e-der2}. Second order terms give condition $p>\frac{Q}{Q-2s}$, while first order terms give $p>\frac{Q}{Q-s}$.
\eproof
\brem Observe that the decomposition given above is not a particular case of results of \cite{moments}. Indeed, there are two main differences with respect to similar results stated in $\R^N$. In the first formula, one can observe that the first-order term $-\sum_{i=1}^n(\int f x_i) (X_i \delta_0)$ contains derivatives in $\L_1$ only, and not in the whole algebra $\L$ as it is in $\R^N$. The second difference is related to the remainder $\sum_{i,j=1}^N X^i X^j(F_{ij})+\sum_{i=n+1}^N X^i(F_i)$, that contains both second order derivatives with respect to all vector fields, and first order derivatives of vector fields outside $\L_1$. Both results can be generalized by defining polynomials on Carnot groups, with the degree being indeed the degree of homogeneity. For details, see \cite{lanconelli-book}.
\erem

\brem Comments given in Remark \ref{rem-RN} apply to this case too. The computation of higher order derivatives is very simple for isotropic dilations, as shown both in \cite{lanconelli-book} and \cite{moments}. But also in this case, we want to highlight the different weights of the vector fields $X^i$.
\erem

\brem
Taylor formula for $\phi$ is indeed a particular case of formulas in \cite[Sec 20.3.2]{lanconelli-book} and \cite{taylor}. We have also simplified the notation, using coordinate $x_i$ instead of exponentiation and logarithm.
\erem

\brem
Following \cite{moments}, one could ask if a decomposition of the kind
\bqn
f=(\int_G f)\delta_0 +\sum_{i,j=1}^N X^i X^j(F_{ij})+\sum_{i=n+1}^N X^i(F_i)
\eqnl{e-2-2}
can be found. We apply $f-(\int f) \delta_0$ to $\phi$ and use \r{e-taylor-kj} with $k=j=1$. We find formula \r{e-2-2} with $$F_{ij}(x):=-w^i w^j x_i x_j \int_0^1 \int_0^1 f\Pt{\Gamma_{\Pt{\frac{r}{t}}^{-1}}(x)} \Pt{\frac{r}{t}}^{-Q-1} \,dr\,\frac{dt}{t}$$ and $$F_i:=w^i (w^i-1) x_i \int_0^1 \int_0^1 f\Pt{\Gamma_{\Pt{\frac{r}{t}}^{-1}}(x)} \Pt{\frac{r}{t}}^{-Q-1} \,dr\,\frac{dt}{t}.$$ We have that $F_{ij}\in L^p(\G)$ when $f(x) x_i x_j\in L^p(\G)$ and $-1<w^i+w^j-1-Q+\frac{Q}{p}<0$. This formula must hold for all $w_i,w_j$. In particular, $p>\frac{Q}{Q-2s+1}$ and $p<\frac{Q}{Q-2}$. This implies that $s<2$, i.e. $s=1$. It means that $G$ is a 1-step Carnot group, i.e. a standard Euclidean space. Thus, the decomposition coincides with the one given in \cite[Thm 2 (a)]{moments} with $j=k=1$. In other terms, this decomposition only holds for Euclidean spaces.
\erem

\section{Large time behavior for the Cauchy problem with hypoelliptic heat diffusion}
\label{s-pde}

We present here the large time behavior of the solutions of \r{e-cauchy}.
\bt
\label{t-nostro1} Let $f_0\in L^1(\G)$ satisfying $|x|_R f\in L^p(\G)$ with $p<\frac{Q}{Q-1}$. Let $f(t)$ be the solution at time $t>0$ of \r{e-cauchy}. Then, given $p\leq q \leq \infty$, there exists $C$, not depending on $t$ and $f_0$, such that
$$\|f(t)- \Pt{\int f_0} P_t\|_q\leq C t^{-Q((1/p)-(1/q))/2-1/2}  \|\, |x|_R f_0\|_p.$$
\et
\bproof
We have $f_0=\Pt{\int f_0} \delta_0 +\mathrm{div} F$, where $F_i\in L^p(\G)$. Recall that $f(t)=f_0\ast P_t$, thus
\bqn
\|f(t)-\Pt{\int f_0} P_t\|_q=\|\sum_{i=1}^N {X^i}F_i \ast P_t\|_q=\|\sum_{i=1}^N F_i \ast {X^i} P_t\|_q\leq \sum_{i=1}^N\|F_i\|_p\|{X^i} P_t\|_r,
\eqnn
where $\frac1p+\frac1r=\frac1q+1$. Recall that $\|F_i\|_p\leq C_{p,\G} \|\, |x|_R f_0\|_p$ and that each $X^i$ is homogeneous of order greater than or equal to 1. Hence $\|f(t)-(\int f_0) P_t\|_q\leq \sum_{i=1}^N C_{p,\G}\|\, \|x\|_R f\|_p t^{-Q/(2r^*)-1/2}\|{X^i} P_1\|_r$, with $\frac1{r^*}=\frac1p-\frac1q$. All the terms $\|{X^i} P_1\|_r$ are finite (Lemma \ref{l-P1}). Consequently, the statement of the theorem holds with $C:=\sum_{i=1}^N C_{p,\G}\|{X^i} P_1\|_r$.
\eproof

We prove now a finer estimate, using Theorem \ref{t:ez2}.
\bt
\label{t-nostro2}
Let $f_0\in L^1 (\G,1+|x|_R)$ such that $|x|_R f,|x|_R^2 f \in L^p(\G)$ with a given $p<\frac{Q}{Q-1}$. Let $f(t)$ be the solution at time $t>0$ of \r{e-cauchy}. Then, given $p\leq q \leq \infty$, there exists $C$, not depending on $t$ and $f_0$, satisfying
\bqn
\|f(t)-\Pt{(\int f_0) P_t+\sum_{i=1}^n (\int f_0 x_i) X^i(P_t)}\|_q\leq C  t^{-Q((1/p)-(1/q))/2-1}  \| |x|^{k+1} f_0\|_p.
\eqnn
\et
\bproof
Using Theorem \ref{t:ez2}, we have 
\bqn
&&\|f(t)-\Pt{(\int f_0) P_t+\sum_{i=1}^n (\int f_0 x_i) X^i(P_t)}\|_q=\nn
&&=\|\sum_{i,j=1}^N X^i X^j(F_{ij})P_t+\sum_{i=n+1}^N X^i(F_i)P_t\|_q\leq\nn
&&\leq \sum_{i,j=1}^N \|F_{ij}\|_p \| X^i X^j (P_t)\|_r+\sum_{i=n+1}^N \|F_i\|_p \|X^i(P_t)\|_r\leq\nn
&&\sum_{i,j=1}^N\| |x|_R^2 f_0\|_p t^{-Q/(2r^*)-1} \|X^iX^j P_1\|_r +\sum_{i=n+1}^N\| |x|_R f_0\|_p t^{-Q/(2r^*)-1} \|X^i P_1\|_r,
\eqnn
where $\frac1p+\frac1r=\frac1q+1$, thus $\frac1{r^*}=\frac1p-\frac1q$. We have used here Lemma \ref{l-P1}, recalling that derivatives $X^iX^j$ are homogeneous of order at least 2, and that the same holds for $X^i$ with $i=n+1,\ldots,N$.
\eproof

This expression shows that, for the large time behavior, the main contribution is given by $P_t$,  its derivatives $X^1 P_t,\ldots, X^n P_t$ in $\L_1$ (case $w^i=1$) and, in general, by terms with small degree of homogeneity. This is a direct consequence of the anisotropicity of the  diffusion on stratified Lie groups.

Higher-order estimations can be easily proved with the same techniques, provided higher-order decomposition of $f_0$ of the kind of Theorem \ref{t:ez2}.\\

\subsection{The Heisenberg group}
\newcommand{\Heis}{{H_2}}
We now apply the results presented before to the \b{Heisenberg group}, defined in Examples \ref{ex-1}, \ref{ex-2}, \ref{ex-3}. The hypoelliptic Laplacian $\Dh$ is the sum of squares $$\Dh\phi=\Pt{\Pt{X^1}^2+ \Pt{X^2}^2}\phi.$$ The explicit expression of $P_t$ has been first computed in \cite{gaveau,hul}. We give here the expression with respect to the definition of $H_2$ given above, see \cite{nostro-kern}:
$$P_t(x,y,z)=\frac{1}{(2\pi t)^2}\int_\R \frac{2 \tau}{\sinh(2\tau)}\exp\left({-\frac{\tau(x^2+y^2)}{2t\tanh(2\tau)}}\right)\cos(2\frac{z\tau}{t})$$

We consider a given function $f_0\in L^1(\Heis,1+|x|_R)$ such that $|x|_R f_0, \|x\|_R^2 f_0\in L^p$ for a fixed $p<\frac43$. We choose $p=1$.  Thus $f_0$ satisfies the hypotheses of Theorem \ref{t:ez2}. We compute the momenta $A_i$ defined as follows
\bqn
A_0&:=&\int_\Heis |f_0(x,y,z)|\, dx\,dy\,dz,\nn
A_1&:=&\int_\Heis |f_0(x,y,z)x|\, dx\,dy\,dz,\quad
A_2:=\int_\Heis |f_0(x,y,z)y|\, dx\,dy\,dz.\quad
\eqnn

We now apply Theorem \ref{t-nostro2} with $q=\infty$, that provides the large time expression in $L^\infty(\Heis)$ for the solution of \r{e-cauchy} on $\Heis$:
$$f(t)=A_0\,P_t+A_1\, X^1(P_t)+A_2\,X^2(P_t)+O(t^{-2/p-1}).$$




{\bf Acknowledgements:} The author is grateful to D. Barilari, C. Mora-Corral and E. Zuazua for useful discussions. He also acknowledges the useful suggestions by the anonymous reviewer.

\end{document}